\documentclass[12pt]{amsart}
\usepackage{amsfonts,amsmath,amssymb}

\advance\textheight by \topskip
\numberwithin{equation}{section}

\title[On Cantor's singular moments]
{On Cantor's singular moments}

\author[H. Prodinger]{Helmut Prodinger}
\address{Helmut Prodinger\\
Centre for Applicable Analysis and Number Theory\\
 Department of Mathematics\\
University of the Witwatersrand, P.~O. Wits\\ 
2050 Johannesburg, South Africa}
\email{helmut@gauss.cam.wits.ac.za}

\date{April 14, 1999}

\begin{document}
\begin{abstract}
We evaluate a constant explicitly, thereby answering a question
raised in \cite{AMM99}.
\end{abstract}

\maketitle

In problem {\bf 10621} of the American Mathematical Monthly, 
Cantor's singular  moments 
$J_n$ were to be computed. In the published
answers \cite{AMM99} they come out as
\begin{equation*}
J_n=\frac{2}{3(n+1)}\sum_{j=0}^n
\binom{n+1}{j}\frac{B_j}{3\cdot2^{j-1}-1}\quad
\text{for}\quad n\ge1
\end{equation*}

and $J_0=1$, with Bernoulli numbers $B_n$.

The editor asked, whether it is possible to compute
\begin{equation*}
J_{-1}=\sum_{n\ge0}J_n
\end{equation*}

exactly.

The purpose of this note is to do that. In \cite{GrPr98} we considered
a similar problem, and the gentle reader
is invited to consult this paper for
more background about the technique, as well as \cite{KnPr96, GrPr96} 
for more information about the Cantor distribution.

Following the method described in \cite{FlSe95}, we can write
$J_n$ as a contour integral viz.
\begin{align*}
J_n&=\frac{2}{3(n+1)}\cdot
\frac{1}{2\pi i}\int_{-\frac12-i\infty}^{-\frac12+i\infty}
\frac{\Gamma(n+2)\Gamma(1-s)}{\Gamma(n+2-s)}
\frac{\zeta(1-s)}{3\cdot2^{s-1}-1}\,ds\\
&=\frac{2}{3}\cdot
\frac{1}{2\pi i}\int_{-\frac12-i\infty}^{-\frac12+i\infty}
\frac{\Gamma(n+1)\Gamma(1-s)}{\Gamma(n+2-s)}
\frac{\zeta(1-s)}{3\cdot2^{s-1}-1}\,ds.
\end{align*}

Therefore
\begin{align*}
\sum_{n=0}^NJ_n&=J_0+
\frac{2}{3}\cdot
\frac{1}{2\pi i}\int_{-\frac12-i\infty}^{-\frac12+i\infty}
\frac{\Gamma(N+2)\Gamma(1-s)}{\Gamma(N+2-s)s}
\frac{\zeta(1-s)}{3\cdot2^{s-1}-1}\,ds\\
&-\frac{2}{3}\cdot
\frac{1}{2\pi i}\int_{-\frac12-i\infty}^{-\frac12+i\infty}
\frac{\Gamma(2)\Gamma(1-s)}{\Gamma(2-s)s}
\frac{\zeta(1-s)}{3\cdot2^{s-1}-1}\,ds.
\end{align*}

From this form, one could even
compute the asympotics as $N\to\infty$. However, here, we only
have to note that the first integral is of order $N^{1-\log_23}$,
which means that it goes to zero. Consequently
\begin{align*}
\sum_{n\ge0}J_n&=1
-\frac{2}{3}\cdot
\frac{1}{2\pi i}\int_{-\frac12-i\infty}^{-\frac12+i\infty}
\frac{1}{(1-s)s}
\frac{\zeta(1-s)}{3\cdot2^{s-1}-1}\,ds\\
&=\frac43
+\frac{2}{3}\cdot
\frac{1}{2\pi i}\int_{\frac32-i\infty}^{\frac32+i\infty}
\frac{1}{s(s-1)}
\frac{\zeta(s)}{3\cdot2^{-s}-1}\,ds\\
&=\frac43
+\frac{2}{3}\sum_{k,m\ge1}3^{-k}\cdot
\frac{1}{2\pi i}\int_{\frac32-i\infty}^{\frac32+i\infty}
\frac{1}{s(s-1)}
\Big(\frac{2^k}{m}\Big)^s\,ds.
\end{align*}
The last step was by using the Dirichlet series for $\zeta(s)$
and the geometric series, both valid for $\Re s=\frac32$.
A simple application of residue calculus,
as it is often used in the context of the Mellin--Perron summation formula
(see \cite{Tenenbaum95, FGKPT94}) evaluates the integrals inside the summation:
\begin{equation*} 
\frac 1{2\pi i}
\int_{\frac32-i\infty}^{\frac32+i\infty}  
\frac1{s(s-1)}t^s ds=
\cases t-1&\text{for }t\geq1\\
0&\text{for }t<1\endcases\;.
\end{equation*}
Therefore
\begin{align*}
\sum_{n\ge0}J_n&=1
+\frac{2}{3}\sum_{k\ge1}
\sum_{1\le m\le 2^k}
3^{-k}
\Big(\frac{2^k}{m}-1\Big)\\
&=1
+\frac{2}{3}\sum_{k\ge1}
\Big(\frac23\Big)^kH_{2^k}-
\frac{2}{3}\sum_{k\ge1}
\Big(\frac23\Big)^k\\
&=-\frac13+\frac{2}{3}\sum_{k\ge1}
\Big(\frac23\Big)^kH_{2^k}=
3.36465\;07281\;00925\;16083\;89349\;6289
\dots,
\end{align*}
with harmonic numbers $H_n=\sum_{1\le k\le n}\frac1k$.

\bibliographystyle{plain}

\end{document}